\documentclass{amsart}
\usepackage{amssymb}
\usepackage{color}

 \DeclareMathOperator{\id}{id}

\copyrightinfo{2006}{American Mathematical Society}

\theoremstyle{definition}

\theoremstyle{remark}

\numberwithin{equation}{section}

\newif \iflv
\long\def\lv #1\endlv{%
    \iflv
    \smallskip
    \begingroup

     \Blue $^{\textstyle **}$\Black~{\footnotesize  #1}~ \Blue $_{\textstyle **}$\Black
    \endgroup
    \smallskip 
    \noindent
    \fi
}

\newif \ifnew

 \lvtrue
\newtrue

\begin{document}

\title{Abelian gradings in Lie algebras}



\author{Esther Garc\'\i a}
\address{ Departamento de Matem\'{a}tica  Aplicada,
Universidad Rey Juan Carlos, 28933 M\'{o}stoles (Madrid), Spain}
\email{esther.garcia@urjc.es}
\thanks{The  first author was partially supported by the MEC and Fondos
FEDER MTM2007-62390 and MTM2010-16153, by FMQ 264, and by MICINN-I3-2010/00075/001.}

\author{Miguel G\'omez Lozano}
\address{Departamento de \'Algebra, Geometr\'{\i}a y
Topolog\'{\i}a, Universidad de M\'alaga, 29071 M\'alaga, Spain}
\email{magomez@agt.cie.uma.es}
\thanks{The second author was partially supported by the MEC and Fondos
FEDER MTM2007-61978 MTM2010-19482, by FMQ 264 and FQM 3737, and by MICINN-I3-2010/00075/001.}

\date{}

\subjclass[2000]{Primary 17B05; Secondary 17B60.}

\keywords{}

\begin{abstract} Given a  Lie algebra $L$ graded by a group $G$, if $L$ is does not contain orthogonal graded ideals and $G$ is generated by the support of $L$, then $G$ is an abelian group.
 \end{abstract}

 \maketitle



It is  well-known  that the subgroup generated by the support of a simple $G$-graded Lie algebra is always abelian. Nevertheless, this fact can be easily extended to Lie algebras without orthogonal graded ideals, and this is the aim of this short note.

Throughout  this note we deal with Lie algebras $L$ over arbitrary rings of scalars $\Phi$ with $\frac 1 2\in \Phi$, with Lie bracket denoted by $[\ ,\ ]$.
Let $G$ be a group. We say that a Lie algebra $L$  is graded
by $G$   if there exists a decomposition $$L=\bigoplus_{g\in G} L_g$$ where each $L_g$ is
a $\Phi$-submodule of $L$ satisfying $[L_g, L_{g'}]\subset L_{gg'}$,  for every
$g,g'\in G$.   The support of an element  $a=\sum_{g\in G} a_g\in L$ is the finite  set ${\rm supp}(a)=\{g\in G\ \vert\
a_g\ne 0\}$, and the support of $L$ as a $G$-graded  algebra  is the set
${\rm supp}(L)=\bigcup_{a\in L} {\rm supp}(a)$.

The next  proposition relates  noncommutative elements of $G$ with orthogonal ideals of  a $G$-graded Lie algebra $L$. We say that a $G$-graded Lie algebra $L$ is graded-prime if it does not contain graded  ideals $I, J\triangleleft L$ such that $[I,J]=0$.

\medskip
\noindent {\bf Proposition} {\it
Let $L$ be a Lie algebra graded by a group $G$ and let $g,g'\in G$ such that $gg'\ne
g'g$. Then $[\id_L(L_g),\id_L(L_{g'})]= 0$. In particular, if $L$ is graded-prime  and
$G$ is generated by the support of $L$, $G$ is an abelian group.}

\begin{proof}
Let us prove the following property: Let $g_1,g_2,\dots, g_n\in G$   $$\hbox{if}\
[L_{g_1},[L_{g_2},[\cdots,[L_{g_{n-1}},L_{g_n}]]]]\ne 0\quad \hbox{then}\
g_ig_j=g_jg_i\ \forall\, i,j\in \{1,2,\dots, n\}.\eqno{(*)}$$ This property is true for
$n=2$:
$$L_{g_1g_2}\supset[L_{g_1},L_{g_2}]=[L_{g_2},L_{g_1}]\subset L_{g_2g_1}$$
which implies, if $0\ne [L_{g_1},L_{g_2}]\subset L_{g_1g_2}\cap L_{g_2g_1}$,  that
$g_1g_2=g_2g_1$. Let us suppose that $(*)$ is true for every $n-1$ elements and let $0\ne
[L_{g_1},[L_{g_2},[\cdots,[L_{g_{n-1}},L_{g_n}]]]]$. By hypothesis,
$$g_ig_j=g_jg_i\ \hbox{if}\ i,j\in \{2,3,\dots,n\}\eqno{(1)}$$ By the Jacobi identity we have
two possibilities:

a).  If  $0\neq [[L_{g_1},L_{g_2}],[\cdots,[L_{g_{n-1}},L_{g_n}]]]\subset
[L_{g_1g_2},[\cdots,[L_{g_{n-1}},L_{g_n}]]]$  then $g_1g_2=g_2g_1$ and, by the
induction hypothesis and (1), the elements $g_1g_2$ and $g_2$ commute with every $g_k$, $k\ge 3$.
Then $g_k(g_1g_2)=(g_1g_2)g_k=g_1g_kg_2$ and multiplying by $g_2^{-1}$ on the right we get $g_1g_k=g_kg_1$ for every  $k\ge 3$.

b). If $[L_{g_2},[L_{g_1},[\cdots,[L_{g_{n-1}},L_{g_n}]]]]\ne 0$, take $h=g_3\dots g_n$. By induction,
$$g_1h=hg_1,\qquad  g_2(g_1h)=(g_1h)g_2, \quad\hbox{ and } g_1g_k=g_k \hbox{ for every $k\ge 3$,}$$
so by (1), $g_2(g_1h)=(g_1h)g_2=g_1g_2h$, which implies $g_1g_2=g_2g_1$.

Now, if $g,g'\in G$ with $gg'\ne g'g$, then for every $g_1,g_2,\dots, g_n\in G$
$$[L_{g'},[L_{g_1},[\cdots[L_{g_n},L_g]]]]=0$$  which proves that  $[\id_L(L_{g'}),\id_L(L_{g})]= 0$.
\end{proof}


\end{document}